\documentclass[10pt]{amsart}

\usepackage{amssymb} 
\usepackage[mathscr]{eucal} 
\usepackage[matrix,arrow,tips,curve,ps]{xy}
\usepackage{amsmath} 
\usepackage{epsfig}
\usepackage{amscd}
\newtheorem{TEO}{Theorem}[section]

\newtheorem{LEM}[TEO]{Lemma}

\newtheorem{COR}[TEO]{Corollary}
\newtheorem{REM}[TEO]{Remark}

\theoremstyle{definition}

\def\OO{{\mathcal O}}

\newcommand\dual{\mathrel{\raise3pt\hbox{$\underline{\mathrm{\thinspace d
\thinspace}}$}}}

\newcommand\proj{\mathbb P}

\newcommand\Co{\mathbb C}

\def\Co{{\mathbb C}}
\begin{document}

\footnote{
Partially supported by  PRIN  2009: "Moduli, strutture geometriche e loro applicazioni" and by INdAM (GNSAGA).
 }

\title{}

\author[E. Colombo]{Elisabetta Colombo}
\address{Dipartimento di Matematica,
Universit\`a di Milano, via Saldini 50,
     I-20133, Milano, Italy } \email{{\tt
elisabetta.colombo@unimi.it}}

\author[P. Frediani]{Paola Frediani}
\address{ Dipartimento di Matematica, Universit\`a di Pavia,
via Ferrata 1, I-27100 Pavia, Italy } \email{{\tt
paola.frediani@unipv.it}}

\title{On the Koszul cohomology of canonical and Prym-canonical binary curves}

\maketitle

\setlength{\parskip}{.1 in}

\begin{abstract}
In this paper we study Koszul cohomology and the Green and Prym-Green
conjectures for canonical and Prym-canonical binary curves.  We prove that
if property $N_p$ holds for a canonical or a Prym-canonical binary curve of genus
$g$ then it holds for a generic canonical or Prym-canonical binary curve of
genus $g+1$. We also verify 
the Green and Prym-Green conjectures for generic canonical and Prym-canonical binary curves of low genus ($6\leq g\leq 15$, $g\neq 8$ for Prym-canonical and $3\leq g\leq 12$ for canonical).

\end{abstract}

\section{Introduction}

Let $C$ be a smooth curve, $L$ a line bundle and $\mathcal{F}$ a
coherent sheaf  on $C$. We recall that the Koszul cohomology group
$K_{p,q}(C,\mathcal{F},L)$ is the middle term cohomology of the
complex:
\begin{equation}
\Lambda^{p+1}H^0(L)\otimes H^0(\mathcal{F}\otimes
L^{q-1})\stackrel{d_{p+1,q-1}}\rightarrow \Lambda^{p}H^0(L)\otimes
H^0(\mathcal{F}\otimes L^{q})\stackrel{d_{p,q}}\rightarrow
\Lambda^{p-1}H^0(L)\otimes H^0(\mathcal{F}\otimes L^{q+1})
\end{equation}
where
$$d_{p,q}(s_1\wedge...\wedge s_p\otimes u):=\sum_{l=1}^{p}(-1)^l s_1\wedge...\wedge \hat{s_l} \wedge...\wedge s_p\otimes
(s_lu).$$ If $\mathcal{F}=\mathcal{O}_C$ the groups
$K_{p,q}(C,\mathcal{O}_C,L)$ are denoted by $K_{p,q}(C,L)$. The
Koszul cohomology theory has been introduced in \cite{gr} and has
been extensively studied in particular in the case of the
canonical bundle. We recall that Green and Lazarsfeld (\cite{gr})
proved that for any smooth curve $C$ of genus $g$ and Clifford
index $c$, $K_{g-c-2,1}(C,K_C)\neq 0$. Green's conjecture says
that this result is sharp i.e. $K_{p,1}(C,K_C)=0$ for all $p\geq
g-c-1$. The Clifford index for a general curve is $[\frac{g-1}{2}]$, so generic Green's
conjecture says that $K_{p,1}(C,K_C)=0$ for all
$p\geq[\frac{g}{2}]$, or equivalently, by duality,
$K_{p,2}(C,K_C)=0$, i.e. property $N_p$ holds, for all
$p\leq[\frac{g-3}{2}]$. Generic Green's conjecture
has been proved by Voisin in \cite{v02},\cite{v05}.
Green's conjecture has also been verified for curves of odd genus and maximal Clifford index (\cite{v05}, \cite{hr}), for general curves of given gonality (\cite{v02}, \cite{te} \cite{sch}), for curves on $K3$-surfaces (\cite{v02}, \cite{v05}, \cite{af}), and in other cases (see \cite{an}).

Another interesting case is when the line bundle is Prym
canonical, $L=K_C\otimes A$ where $A$ is a non trivial 2 torsion
line bundle. This case has been studied in \cite{fl}, where the
Prym-Green conjecture has been stated. This is an analogue of the
Green conjecture for general curves, namely it says that for a
general Prym-canonical curve $(C,K_C\otimes A)$, we have
$K_{p,2}(C,K_C\otimes A)=0$, i.e. property $N_p$ holds, for all
$p\leq[\frac{g}{2}-3]$. Prop.3.1 of \cite{fl} shows that for any
$(C,K_C\otimes A)$ and $p>[\frac{g}{2}-3]$, $K_{p,2}(C,K_C\otimes
A) \neq 0$.

Debarre in \cite{deba} proved that a generic Prym-canonical curve
of genus $g \geq 6$ is projectively normal (property $N_0$) and
for $g \geq 9$ its ideal is generated by quadrics (property
$N_1$).

In \cite{cfes} the Prym-Green conjecture is proved for genus
$g=10,12,14$ by degeneration to irreducible nodal curves and
computation with Macaulay2. In a private communication Gavril Farkas told
us that they could verify the conjecture also for $g=18,20$. The
computations made in \cite{cfes} for genus 8 and 16 suggest that the
Prym-Green conjecture may be false for genus which is a multiple of 8 or
perhaps a power of 2. The possible failure of the Prym-Green
conjecture in genus 8 is extensively discussed in the last section
of \cite{cfes}, where a geometric interpretation of this phenomenon is given.

In this paper we study Koszul cohomology and the Green and Prym-Green
conjectures for canonical and Prym-canonical binary curves. Recall
that a binary curve of genus $g$ is a stable curve consisting of
two rational components $C_j$, $j=1,2$ meeting transversally at
$g+1$ points. The canonical and Prym-canonical models of binary
curves that we analyze are the one used in \cite{ccm} and
\cite{cf} and described in the next section. The
main result of the paper (Theorem \eqref{indstep}) says that
if property $N_p$ holds for a Prym-canonical binary curve of genus
$g$ then it holds for a generic Prym-canonical binary curve of
genus $g+1$. In particular,  if the Prym-Green conjecture is true for a Prym-canonical binary curve of genus $g = 2k$, then it is true for a general
Prym-canonical binary curve of genus $g = 2k+1$.

Moreover we verify
the conjecture by a direct computation for $g=6,9,10,12,14$ (see Corollary \eqref{corpg}).

As a consequence, we show that the generic Prym-canonical curve of
genus $g$ satisfies property $N_0$ for $g \geq 6$,  property $N_1$
for $g \geq 9$ (already shown by Debarre),  property $N_2$ for $g
\geq 10$, property  $N_3$ for $g \geq 12$ and property $N_4$ for
$g \geq 14$ (Corollary \eqref{Np}).

For $g =8$ and $g =16$ our computations on Prym-canonical binary curves also suggest that Prym-Green conjecture's might fail, in fact in our examples we find that $K_{\frac{g}{2}-3,2}(C, K_C \otimes A) =1$ both for $g=8$ and $g=16$ (see Remark \eqref{g8-16}).

An analogous result of Theorem \eqref{indstep} is proven for canonically embedded binary curves (Theorem \eqref{can}), where we show that  if property $N_{p}$ holds for a canonical binary curve of genus $g$, then the same property holds for a general canonical binary curve of genus $g+1$. In particular,  if the Green conjecture is true for a canonical binary curve of genus $g = 2k-1$, then it is true for a general
canonical binary curve of genus $g = 2k$.

Theorem \eqref{indstep} and analogous computations with maple in genus $g =3,5,7, 9,11$, imply that for a general canonical binary curve, if  $g \geq 3$, then propery $N_0$ holds (see also \cite{ccm} section 2), if $g \geq 5$, then propery $N_1$ holds, if $g \geq 7$, then propery $N_2$ holds, if $g \geq 9$, then property $N_3$ holds, and if $g \geq 11$, then property $N_4$ holds.

{\bf Acknowledgments.} We thank Riccardo Murri for having been so kind to do for us the computer computations in $g=14,16$.

\section{Canonical and Prym-canonical binary curves}
\subsection{Construction of canonical binary curves}

Recall that a binary curve of genus $g$ is a stable curve
consisting of two rational components $C_j$, $j=1,2$ meeting
transversally at $g+1$ points. Moreover, $H^0(C,\omega_C)$ has dimension $g$ and the restriction of
$\omega_C$ to  the component $C_j$ is $K_{C_j}(D_j)$
where $D_j$ is the divisor of nodes on $C_j$. Since
$K_{C_j}(D_j)\cong \OO_{{\proj}^{1}}(g-1)$ we observe that the
components are embedded by the complete linear system
$|\OO_{{\proj}^{1}}(g-1)|$ in ${\proj}^{g-1}$.

Following \cite{ccm}, we assume that the first $g$ nodes are
$P_i=(0,...,0,1,0,...0),$ with 1 at the $i$-th place, $i=1,...,g$.  Then we can
assume that $C_j$ is the image of the map
\begin{equation}\label{can}
\begin{gathered}\phi_j:{\proj}^1 \rightarrow {\proj}^{g-1}, \
j=1,2\\
\phi_j(t,u):= [\frac{M_j(t,u)}{(t-a_{1,j}u)},
...,\frac{ M_j(t,u)}{(t-a_{g-1,j}u)}]
\end{gathered}
\end{equation}
with $M_j(t,u):= \prod_{r=1}^{g} (t-a_{r,j}u)$, $j=1,2$ and
$\phi_j([a_{l,j},1]) = P_l$, $l=1,...,g$.

We  see  that the remaining node is the point
 $P_{g+1}:=[1,...,1]$ and it is the image of $[1,0]$  through the maps
 $\phi_j$, $j=1,2$. One can easily check that, for generic values of the $a_{i,j}$'s, $C=C_1\cup
 C_2$ is a canonically embedded binary curve.

\subsection{Construction of Prym-canonical binary curves}

Let $C$ be a binary curve of genus $g$, and $A\in Pic^0(C)$ a nontrivial line bundle. Then $H^0(C,\omega_C \otimes A)$ has dimension $g-1$
and  the restriction of $\omega_C \otimes A$ to  the component $C_j$ is
$K_{C_j}(D_j)$ where $D_j$ is the divisor of nodes on $C_j$. Since
$K_{C_j}(D_j)\cong \OO_{{\proj}^{1}}(g-1)$, the
components are embedded by a linear subsystem of
$\OO_{{\proj}^{1}}(g-1)$, hence they are projections from a point
of rational normal curves in ${\proj}^{g-1}$. Viceversa, let us
take 2 rational curves embedded in  ${\proj}^{g-2}$ by non
complete linear systems of degree $g-1$ intersecting transversally
at $g+1$ points. Then their union $C$ is a binary curve of genus
$g$ embedded either by a linear subsystem of $\omega_C$ or by a
complete linear system $|\omega_C \otimes A|$, where $A\in Pic^0(C)$
is nontrivial (see e.g. \cite{capo}, Lemma 10). In \cite{cf} (Lemma 3.1) we  constructed a
binary curve $C$ embedded in ${\proj}^{g-2}$ by  a linear system
$|\omega_C \otimes A|$ with $A^{\otimes 2}\cong \OO_C$, and $A$ is non
trivial. Let us now recall this construction and denote a binary curve with this embedding a Prym-canonical binary curve.

Assume that the first $g-1$ nodes, are $P_i=(0,...,0,1,0,...0)$ with 1 at the $i$-th place, $i=1,...,g-1$, the remaining two nodes are $P_g:=[t_1,...,t_{g-1}]$ with $t_i=0$ for $i
=1,...,[\frac{g}{2}]$, $t_i =1$, for $i =
[\frac{g}{2}]+1,...,g-1$.
and $P_{g+1}:=[s_1,...,s_{g-1}]$ with $s_i=1$ for $i
=1,...,[\frac{g}{2}]$, $s_i =0$, for $i =
[\frac{g}{2}]+1,...,g-1$.

Then the component $C_j$ is the image of the map

\begin{equation}\label{pcan}
\begin{gathered}\phi_j:{\proj}^1 \rightarrow {\proj}^{g-2}, \
j=1,2, \ \text{where} \\
\phi_1(t,u):= [\frac{tM_1(t,u)}{(t-a_{1,1}u)},..., \frac{tM_1(t,u)}{(t-a_{k,1}u)}, \frac{-M_1(t,u)d_1 a_{k+1,1}u}{A_1(t-a_{k+1,1}u)},..., \frac{-M_1(t,u)d_1a_{g-1,1}u}{A_1(t-a_{g-1,1}u)}]\\
\phi_2(t,u):= [\frac{tM_2(t,u)}{(t-a_{1,2}u)},..., \frac{tM_2(t,u)}{(t-a_{k,2}u)}, \frac{-M_2(t,u)d_2 a_{k+1,2}u}{A_2(t-a_{k+1,2}u)},..., \frac{-M_2(t,u)d_2a_{g-1,2}u}{A_2(t-a_{g-1,2}u)}]\\
\end{gathered}
\end{equation}
with $k :=[\frac{g}{2}]$, $M_j(t,u):= \prod_{r=1}^{g-1} (t-a_{r,j}u)$,  and
$A_j= \prod_{i=1}^{g-1} a_{i,j}$, $j=1,2$, $d_2$ is a nonzero constant and $d_1 = \frac{-d_2 A_1}{A_2}$.
Notice that we have $\phi_j([a_{l,j},1]) = P_l$, $l=1,...,g-1$, $\phi_j([0,1]) = P_g$, $\phi_j([1,0]) = P_{g+1}$, $j=1,2$.
In Lemma 3.1 of \cite{cf} we proved that for a general choice of $a_{i,j}$'s, $C=C_1\cup
C_2$ is a binary curve embedded in ${\proj}^{g-2}$ by  a linear
system $|\omega_C \otimes A|$ with $A^{\otimes 2}\cong \OO_C$ and $A$ nontrivial. In fact, recall
that $Pic^0(C) \cong {{\Co}^*}^g \cong {{\Co}^*}^{g+1}/{\Co}^*$, where ${\Co}^*$ acts diagonally,
and in Lemma 3.1 of \cite{cf} it is shown and our line bundle $A$ corresponds to the element $[(h_1,...,h_{g+1})] \in {{\Co}^*}^{g+1}/{\Co}^*$,
where $h_i=1$, for $i< [\frac{g}{2}]+1$, $h_i = -1$, for $i
=[\frac{g}{2}]+1,...,g-1$, $h_g=-1$, $h_{g+1} = 1$, so in particular $A$ is of
2-torsion.

\section{Property $N_p$ for Prym-canonical binary curves}
Let $C \subset {\proj}^{g-2}$ be a Prym-canonical binary curve
embedded by $\omega_C \otimes A$, with $A^{\otimes 2} \cong
\OO_C$, as in (\ref{pcan}). In this section we study the Koszul
cohomology for these curves, in particular we investigate property
$N_p$, i.e. the vanishing of $K_{p,2}(C, K_C \otimes A)$. Since by
duality (\cite{gr}, see also \cite{fr} prop.1.4) we have $K_{p,2}(C, K_C \otimes A)
\cong K_{g-3-p,0}(C, K_C,K_C \otimes A)^{\vee} $,  this vanishing
is equivalent to the injectivity of the Koszul map
\begin{equation}
\label{koszul}
F_{g-3-p}: \Lambda^{g-3-p}H^0(C,\omega_C \otimes A) \otimes H^0(C,\omega_C) \rightarrow \Lambda^{g-4-p}H^0(C,\omega_C \otimes A) \otimes H^0(C,\omega_C^2 \otimes A).
\end{equation}

Our strategy is to compare this map with analogous Koszul maps for a partial normalization of the curve $C$ at one node and possibly use induction on the genus.

To this end, let us introduce some notation: set $k :=[\frac{g}{2}]$ and denote by $\tilde{C}_r$ the partial normalization of $C$ at the node $P_r$ with $r\leq k$ if $g =2k$,  $r \geq k+1$ if $g = 2k+1$. This choice of the node is necessary in order to obtain the Prym-canonical model for the curve $\tilde{C}_r$. In fact, observe that in this way, for a general choice of the $a_{i,j}$'s, the
projection  from $P_r$ sends the curve $C$ to the Prym-canonical
model of $\tilde{C}_r$ in ${\proj}^{g-3}$ given by the line bundle
$K_{\tilde{C}_r} \otimes  A'_r$ where $A'_r$
corresponds to the point $(h'_1,...,h'_{g-1}, 1) \in {{\Co}^*}^{g}/{\Co}^*$, with $h'_i = 1$
for $i\leq [\frac{g-1}{2}]$, $h'_i = -1$ for
$i=[\frac{g-1}{2}]+1,...,g-1$,  as described above. In
fact  $(\tilde{C}_r,A'_r)$ is parametrized by
$a'_{i,j} = a_{i,j}$ for $i
\leq r-1$, $j =1,2$, $a'_{i,j} = a_{i+1,j}$ for $i
\geq r$, $j =1,2.$ So if we set $d'_{j} := \frac{d_j}{a_{r,j}}$, $j =1,2$, we clearly have a pair $(\tilde{C}_r,A'_r)$ as in \eqref{pcan}. For simplicity let us choose $d_2 =1$, so $d_1 = -\frac{A_1}{A_2}$, hence $d'_{2} := \frac{1}{a_{r,2}}$, $d'_{1} := -\frac{A_1}{A_2 a_{r,1}}$.

To simplify the notation, set $T_g:= H^0(C,\omega_C \otimes A)$, $H_g :=  H^0(C,\omega_C)$, $B_g:= H^0(C,\omega_C^2 \otimes A)$.
 Denote by $\{t_1,...,t_{g-1}\}$ the basis of $T_g$ given by the coordinate hyperplane sections  in $\proj^{g-2} \cong \proj (T_g^{\vee})$ and by $\{s_1,...,s_{g}\}$ the basis of $H_g$ given by the coordinate hyperplane sections  in $\proj^{g-1} \cong \proj (H_g^{\vee})$.
$T_{g-1,r}:= H^0(\tilde{C}_r,\omega_{\tilde{C}_r} \otimes A'_r)$, $H_{g-1,r} :=  H^0(\tilde{C}_r,\omega_{\tilde{C}_r})$, $B_{g-1,r}:= H^0(\tilde{C}_r,\omega_{\tilde{C}_r}^2 \otimes A'_r)$.
Denote by  $\{t'_1,...,t'_{g-2}\}$ the basis of $T_{g-1,r}$ given by the coordinate hyperplane sections  in $\proj^{g-3} \cong \proj (T_{g-1,r}^{\vee})$ and by $\{s'_1,...,s'_{g-1}\}$ the basis of $H_{g-1,r}$ given by the coordinate hyperplane sections  in $\proj^{g-2} \cong \proj (H_{g-1,r}^{\vee})$.

We have the following injections:
\begin{equation}
T_{g-1,r} \stackrel{I_r}\hookrightarrow T_g, \ t'_i \mapsto t_i \  \text{for} \  i \leq r-1, \ t'_i \mapsto t_{i+1} \  \text{for} \  i \geq r,\end{equation}
\begin{equation}
\label{h}
H_{g-1,r} \stackrel{J_r}\hookrightarrow H_g, \ s'_i \mapsto s_i \  \text{for} \  i \leq r-1,  \ s'_{i}\mapsto s_{i+1} \  \text{for} \  i \geq r.\end{equation}

Clearly these maps induce an injective map
\begin{equation}
B_{g-1,r} \stackrel{L_r}\hookrightarrow B_g,
\end{equation}
which on the set of generators of $B_{g-1,r}$ given by $t'_i s'_j$, $i = 1,...,g-2$, $j = 1,...,g-1$ is given by $t'_i s'_j \mapsto I_r(t'_i)J_r(s'_j)$. We claim that this map is well defined and injective by the definition of the $t'_i$'s and $s'_j$'s. In fact the restriction of $\sum \alpha_{i,j} t'_i s'_j$ to the two rational components of $\tilde{C}_r$ yields two polynomials $Q_1$ and $Q_2$. On the other hand we have $(\sum \alpha_{i,j} I_r(t'_i) J_r(s'_j))_{|C_{i}} = (t-a_{r,i})^2 Q_i$, hence  $L_r$ is well defined and injective.
We finally have a map
\begin{equation}
\Lambda^{l-1} T_{g-1,r} \stackrel{\wedge t_{r}} \longrightarrow \Lambda^{l} T_g, \end{equation}
where by $\wedge t_{r}$ we indicate the composition of the natural  map induced by $I_r$ at the level of the $(l-1)$-th exterior power $\Lambda^{l-1} T_{g-1,r} \rightarrow \Lambda^{l-1} T_g$ composed by the wedge product  with $t_{r}$, $ \Lambda^{l-1} T_g \stackrel{ \wedge t_{r}}  \longrightarrow \Lambda^{l} T_g $.

As in \eqref{koszul}, denote  by $F_l : \Lambda^l T_g \otimes H_g \rightarrow \Lambda^{l-1} T_g \otimes B_g$ the Koszul map.

We have the following commutative diagram

\begin{equation}
\label{diagram1} \xymatrix{
  \Lambda^{l} T_{g} \otimes H_{g} \ar[r]^{F_{l}}  \ar[r]& \Lambda^{l-1} T_{g} \otimes B_{g} \ar[r]^{  \pi_{r}}
& \langle t_{r} \rangle   \wedge \Lambda^{l-2} T_g \otimes B_g& \\
\Lambda^{l-1} T_{g-1,r} \otimes H_{g-1,r} \ar[r]^{\tilde{F}_{l-1}}
\ar[u]^{\wedge t_{r} \otimes J_r}&  \Lambda^{l-2} T_{g-1,r}
\otimes B_{g-1,r}  \ar[ur]^{\wedge t_{r} \otimes L_r}}
\end{equation}

From now on, given a multi-index $I=(i_1,...,i_l)$ we denote by $t_I:=t_{i_1}\wedge ... \wedge t_{i_l}$.

To study the injectivity of the maps $F_l$, a preliminary reduction comes from the following

\begin{LEM}
\label{W} Let $W  \subset \Lambda^l T_g \otimes H_g$ be the
subspace generated by the elements of the form $t_I \otimes s_j$,
where $j \not \in I$. Then the kernel of the Koszul map $F_l : \Lambda^l T_g \otimes H_g
\rightarrow \Lambda^{l-1} T_g \otimes B_g$ is contained in $W$.
\end{LEM}
\proof
Assume that $ v \in \Lambda^l T_g \otimes H_g $, $ v = \sum_{I, |I|=l} \sum_{j=1...g} \lambda^I_j t_I \otimes s_j$ is such that $F_{l}(v) = 0$.
$F_l(v) = \sum_{J, |J| = l-1}  \sum_{ I =J \cup \{m\}} \sum_{j=1...g} \lambda^I_j \epsilon(I,J) t_J \otimes t_m s_j=0$, where $\epsilon(I,J) = \pm 1$, depending on the position of $m$ in the multi-index $I = J \cup \{m\}$. Then if we fix a multi-index $J$ with $|J| = l-1$, we must have $ \sum_{ m} \sum_{j=1...g} \lambda^{J \cup \{m\}}_j \epsilon(J \cup \{m\},J) t_J \otimes t_m s_j =0$ and therefore $$ \sigma_J:= \sum_{ m} \sum_{j=1...g} \lambda^{J \cup \{m\}}_j \epsilon(J \cup \{m\},J)  t_m s_j =0.$$ So we have ${\sigma_J}_{|C_1} \equiv 0$, namely, if we denote by $P_1(t) := t \cdot M_1(t,1)$, as in \eqref{pcan}, we have
$$\sum_{j=1...g}  \sum_{ m \leq k} \lambda^{J \cup \{m\}}_j \epsilon(J \cup \{m\},J) \frac{ P_1(t)}{t-a_{m,1}} \frac{ P_1(t)}{t-a_{j,1}} $$
$$ +\sum_{j=1...g}  \sum_{ m \geq k+1} \lambda^{J \cup \{m\}}_j \epsilon(J \cup \{m\},J) \frac{ P_1(t) a_{m,1}}{A_2t(t-a_{m,1})} \frac{ P_1(t)}{t-a_{j,1}} =0$$
 If we evaluate in $t = a_{m,1}$, there remains only one term in the sum, namely the one with $j =m$, and hence we have

 $$ \lambda^{J \cup \{m\}}_m  \epsilon(J \cup \{m\},J) a^2_{m,1} \cdot  \prod_{ r \neq m, r =1...g-1} (a_{m,1} - a_{r,1})^2 =0, \ \text{if } \  m \leq k, $$
  $$ \lambda^{J \cup \{m\}}_m  \epsilon(J \cup \{m\},J) \frac{a^2_{m,1}}{A_2} \cdot  \prod_{ r \neq m, r =1...g-1} (a_{m,1} - a_{r,1})^2 =0, \ \text{if } \  m \geq k+1, $$
  hence we have $ \lambda^{J \cup \{m\}}_m =0$ for all $m$.

Since this holds for every multi-index $J$ of cardinality $l-1$, we have shown that we can write $v = \sum_{I, |I|=l} \sum_{j=1...g, j \not \in I} \lambda^I_j t_I \otimes s_j$.
\qed

We can now state and prove our main result.

\begin{TEO}
\label{indstep}
Assume that $g = 2k$, or $g = 2k+1$ and take an integer $p \leq k-3$.
If property $N_p$ holds for a binary curve $\tilde{C}$ of genus $g-1$
embedded in ${\proj}^{g-3}$ by $|\omega_{\tilde{C}} \otimes A'|$ as in \eqref{pcan}
for a generic choice of the parameters $a'_{i,j}$, then it holds for all binary curves $C$
of genus $g$ embedded in ${\proj}^{g-2}$ by $|\omega_C \otimes A|$ as in \eqref{pcan}
for a generic choice of the $a_{i,j}$'s.

\end{TEO}

\proof
We want to prove that $K_{p,2}(C, K_C \otimes A) =0$ for a binary curve of genus $g$ and we know that $K_{p,2}(\tilde{C_r}, K_{\tilde C_r} \otimes A'_r) =0$, for the curve $\tilde{C_r}$ which is obtained from $C$ by projection from  $P_{r}$ with $r \geq k+1$ if $g = 2k+1$,  $r \leq k$ if $g = 2k$.

By duality, $K_{p,2}(C, K_C \otimes A) \cong K_{g-3-p,0}(C, K_C,K_C \otimes A)^{\vee} $, so the statement is equivalent to prove injectivity of the Koszul map
$$ F_{g-3-p}: \Lambda^{g-3-p}T_g \otimes H_g \rightarrow \Lambda^{g-4-p}T_g \otimes B_g.$$
By assumption we know injectivity of the map
$$\tilde{F}_{g-4-p}: \Lambda^{g-4-p}T_{g-1,r} \otimes H_{g-1,r} \rightarrow \Lambda^{g-5-p}T_{g-1,r} \otimes B_{g-1,r}.$$
For simplicity let us denote by $l:=g-3-p$.

Assume first of all that $g = 2k+1$ and consider the projection of
$C$ from $P_{g-1}$.

Recall that by Lemma \eqref{W} we can reduce to prove injectivity of $F_{l}$ restricted the subspace $W$ generated by such $T_I \otimes s_j$ with $j \not \in I$. Note that we can decompose $W$ as  $W:=X_{g-1} \oplus Y_{g-1}$, where $X_{g-1}$ is the intersection with $W$ of the image of the map $\wedge t_{g-1} \otimes J_{g-1}$ in diagram \eqref{diagram1} and $Y_{g-1}$ is the subspace of $W$ generated by such $t_I \otimes s_j$ with $g-1 \not \in I$ and $j \not \in I$: $$X_{g-1} = \langle t_{g-1} \wedge  t_{J}  \otimes s_{j} \ | \   j \not \in J \rangle,  \ Y_{g-1} = \langle t_I \otimes s_{j} \ | \ j ,g-1 \not \in I \rangle.$$

Assume now that $F_{l}(x_{g-1} + y_{g-1}) =0$, where $x_{g-1} \in X_{g-1}$, $y_{g-1} \in Y_{g-1}$.  Then  we have $0=\pi_{g-1} \circ F_{l}(x_{g-1} + y_{g-1}) = \pi_{g-1} \circ F_{l}(x_{g-1})= (\wedge t_{g-1} \otimes L_{g-1}) \circ  \tilde{ F}_{l-1}(x_{g-1})$, by the commutativity of  diagram \eqref{diagram1}. Hence $x_{g-1}=0$, since by induction we are assuming that $\tilde{ F}_{l-1}$ is injective.
 So we have  reduced to prove injectivity of $F_{l}$ restricted to $Y_{g-1}$.

 Now consider the projection of $C$ from the point $P_{g-2}$.

Set
$$Y'_{g-2} = \langle t'_{J}  \otimes s'_{j} \ | \    j, g-2 \not \in J \rangle \subset  \Lambda^{l-1} T_{g-1,g-2} \otimes H_{g-1,g-2}$$

Observe that the image $X_{g-2}:=  (\wedge t_{g-2} \otimes J_{g-2}) (Y'_{g-2})$ is contained in $Y_{g-1}$ and in fact
$$X_{g-2}=  \langle t_{g-2} \wedge  t_{J}  \otimes s_{j} \ | \   j,g-1 \not \in J \rangle.$$
So we have $Y_{g-1} = X_{g-2} \oplus Y_{g-2}$, where $Y_{g-2}$ is the subspace of  $Y_{g-1}$ generated by those elements of the form  $t_I \otimes s_j$ where $g-2,g-1,j \not \in I$.
We have the following commutative diagram

\begin{equation}
\label{diagram2} \xymatrix{
  \Lambda^{l} Y_{g-1} \ar[r]^{F_{l}}  \ar[r]& \Lambda^{l-1} T_{g} \otimes B_{g} \ar[r]^{  \pi_{g-2}}
& \langle t_{g-2} \rangle   \wedge \Lambda^{l-2} T_g \otimes B_g& \\
{ \ \ \ \ \ \ \ \Lambda^{l-1}Y'_{g-2} \ \ \ \ \ \ \ }  \ar[r]^{\tilde{F}_{l-1}}  \ar[u]^{\wedge t_{g-2} \otimes J_{g-2}}&  \ \ \ \ \Lambda^{l-2} T_{g-1,g-2} \otimes B_{g-1,g-2}  \ar[ur]^{\wedge t_{g-2} \otimes L_{g-2}}}
\end{equation}

Assume that $v=x_{g-2} + y_{g-2}   \in Y_{g-1} = X_{g-2} \oplus Y_{g-2}$ is such that $F_{l}(v) =0$, then  we have $0 = \pi_{g-2} \circ F_{l}(x_{g-2}+y_{g-2}) = \pi_{g-2} \circ F_{l}(x_{g-2})$. So $0= \tilde{F}_{l-1}(x_{g-2})$ by the commutativity of the diagram, and this implies $x_{g-2} = 0$ by induction.
Therefore  we can assume that $v \in Y_{g-2}$, hence $v$ is a linear combination of vectors of the form $t_I \otimes s_j$ where $g-2,g-1,j \not \in I$.

Repeat the procedure, i.e. project from the points $P_r$, $r=g-3...l$. This can be done since $l = g-3-p \geq k+1$.  In this way we can reduce to prove injectivity for the  restriction of the map $F_{l}$ to the subspace $Y_{l}$ of $W$ generated by the elements of the form $t_I \otimes  s_j$ where $l,...,g-1,j \not \in I$. Observe that since $|I| = l$, we have  $Y_{l} = 0$, so $F_{l}$ is injective and the theorem is proved.

If $g =2k$ the proof is analogous: we subsequently project from the points $P_1, P_{2},...,P_{g-l}$. As before note that this can be done since $g-l= p+3 \leq k$.  In this way we reduce to prove injectivity for the  restriction of the map $F_{l}$ to the subspace $Y$ of $W$ generated by the elements of the form $t_I \otimes  s_j$ where $1,2,3,...,g-l, j \not \in I$ and since $|I| = l$, we have  $Y= 0$, so $F_{l}$ is injective and the theorem is proved.
\qed

\begin{COR}
If the Prym-Green conjecture is true for a Prym-canonical binary curve of genus $g = 2k$  as in \eqref{pcan}, then it is true for a Prym-canonical binary curve of genus $g = 2k+1$ as in \eqref{pcan} for generic parameters $a_{i,j}$.
\end{COR}
\proof
The conjecture for $g = 2k+1$ says that $K_{k+1,0}(C,K_C,K_C \otimes A) = 0$, or analogously that property $N_{k-3}$ holds for a generic $C$ embedded with $K_C \otimes A$. Hence the corollary immediately follows from Theorem  \eqref{indstep} with $i =k-3$.
\qed

\begin{COR}
\label{corpg}
\label{Np}The generic Prym-canonical curve of genus $g$ satisfies
property $N_0$ for $g \geq 6$, $N_1$ for $g \geq 9$, $N_2$ for $g
\geq 10$, $N_3$ for $g \geq 12$, $N_4$ for $g \geq 14$.
\end{COR}

\proof With a direct computation one verifies the Prym-Green conjecture for
explicit examples of Prym-canonical binary curves as in \eqref{pcan} for $g = 6,
9,10,12,14$, so the proof follows from Theorem  \eqref{indstep}
for generic Prym-canonical binary curves, and then by
semicontinuity for generic Prym-canonical smooth curves.

To do the computations we wrote a very simple maple code ( http://www-dimat.unipv.it/~frediani/prym-can) in which we
explicitly give the matrix representing the Koszul map $F_l$: for
every multi-index $J$ with $|J| = l-1$, we take the projection of
the image of $F_l$ onto $t_J \otimes B_{g}$ and we restrict it to
the rational components $C_j$. So we have two polynomials in one
variable and we take their coefficients.

Once the matrix is constructed, for $g = 6,
9,10,12$, maple computed its rank modulo $131$, which turned out to be maximal.  In the case $g =14$ the order of the matrices was too big, so Riccardo Murri made the rank computation using the Linbox (\cite{1}) and Rheinfall (\cite{2}) free software libraries. Two different rank
computation algorithms were used: Linbox' "black box" implementation
of the block Wiedemann method (\cite{4,5}), and Rheinfall's Gaussian
Elimination code(\cite{6}).  Results obtained by either method agree.

In both cases, the GNU GMP library (\cite{3}) provided the underlying
arbitrary-precision
representation of rational numbers and exact arithmetic operations.

\qed

\begin{REM}
\label{g8-16}
For Prym-canonical curves of genus 8, the maple computation on specific examples of binary curves gives $dimK_{1,2}(C, K_C \otimes A)=1$. This result is compatible with the computations in \cite{cfes}.

For Prym-canonical  binary  curves of genus 16, we constructed the matrix representing the Koszul map $F_8$ on examples using maple and Riccardo Murri computed its rank as explained in the proof of Corollary \eqref{corpg}.  Again it turned out that $dimK_{5,2}(C, K_C \otimes A)=1$, confirming the computations in \cite{cfes}.

\end{REM}

\section{Property $N_p$ for canonical binary curves}

In analogy with the Prym-canonical case, we study  now property
$N_p$ for canonical binary curves with the same inductive method,
projecting from a node. So, let $C \subset {\proj}^{g-1}$ be a
canonical binary curve and denote by $\tilde{C}_r$ the partial
normalization of $C$ at the node $P_r$,   $1\leq r \leq g$. As
above, for a general choice of the $a_{i,j}$'s, the projection
from $P_r$ sends the curve $C$ to the canonical model of
$\tilde{C}_r$ in ${\proj}^{g-2}$, where  $\tilde{C}_r$ is
parametrized by $a'_{i,j} = a_{i,j}$ for $i \leq r-1$, $j =1,2$,
$a'_{i,j} = a_{i+1,j}$ for $i \geq r$, $j =1,2$.

Set $H_g :=  H^0(C,\omega_C)$, $D_g:= H^0(C,\omega_C^2)$, $F_l : \Lambda^l H_g
\otimes H_g \rightarrow \Lambda^{l-1} H_g\otimes D_g$ the Koszul
map. Denote as before by $\{s_1,...,s_{g}\}$ the basis of $H_g$ given by the
coordinate hyperplane sections  in $\proj^{g-1} \cong \proj
(H_g^{\vee})$ .

$H_{g-1,r} :=  H^0(\tilde{C}_r,\omega_{\tilde{C}_r})$, $D_{g-1,r}:= H^0(\tilde{C}_r,\omega_{\tilde{C}_r}^2)$.
Denote by   $\{s'_1,...,s'_{g-1}\}$ the basis of $H_{g-1,r}$ given by the coordinate hyperplane sections  in $\proj^{g-2} \cong \proj (H_{g-1,r}^{\vee})$.

We have the injections $H_{g-1,r} \stackrel{J_r}\hookrightarrow H_g,$ as in \eqref{h} and
$D_{g-1,r} \stackrel{L_r}\hookrightarrow D_g,$
which on the set of generators of $B_{g-1,r}$ given by $s'_i s'_j$,
$i,j = 1...g-1$,  is given by $s'_i s'_j \mapsto J_r(s'_i)J_r(s'_j)$.

We
finally have a map
\begin{equation}
\Lambda^{l-1} H_{g-1,r} \stackrel{\wedge s_{r}} \longrightarrow
\Lambda^{l} H_g, \end{equation} where by $\wedge s_{r}$ we
indicate the composition of the natural  map induced by $J$ at the
level of the $l-1$-th exterior power $\Lambda^{l-1} H_{g-1,r}
\rightarrow \Lambda^{l-1} H_g$ composed by the wedge product  with
$s_{r}$, $ \Lambda^{l-1} H_g \stackrel{ \wedge s_{r}}
\longrightarrow \Lambda^{l} H_g $.

We are interested in property $N_p$ for these curves, hence by duality, in the vanishing of $K_{g-2-p,1}(C,K_C)$.  Clearly the vanishing of $K_{l,1}(C,K_C)$ is equivalent to the injectivity of the map
\begin{equation}
\frac{\Lambda^l H_g \otimes H_g}{ \Lambda^{l+1} H_g}\rightarrow
\Lambda^{l-1} H_g\otimes D_g
\end{equation}
coming from the Koszul complex.

Notice that there is an isomorphism between $\frac{\Lambda^l H_g
\otimes H_g}{ \Lambda^{l+1} H_g}$ and the subspace $V_{g}$ of
$\Lambda^l H_g \otimes H_g$ generated by the elements of the form
$s_I \otimes s_j$, where $j \geq i_1$, so the above injectivity is
equivalent to the injectivity of the restriction of $F_l$ to
$V_{g}$.

We have the following commutative diagram

\begin{equation}
\label{diagram1can} \xymatrix{
  V_g \ar[r]^{F_{l} }  \ar[r]& \Lambda^{l-1} H_{g} \otimes D_{g} \ar[r]^{  \pi_{r}}
& \langle s_{r} \rangle   \wedge \Lambda^{l-2} H_g \otimes D_g& \\
\ \ \ \ \ \ \ V_{g-1,r} \ \ \ \  \ar[r]^{  \tilde{F}_{l-1}} \ar[u]^{\wedge s_{r} \otimes
J_r}& \ \ \ \  \Lambda^{l-2} H_{g-1,r} \otimes D_{g-1,r}  \ar[ur]^{\wedge
s_{r} \otimes L_r}}
\end{equation}

where $V_{g-1,r}$ is the subspace of $\Lambda^{l-1} H_{g-1,r} \otimes  H_{g-1,r}$ generated by the elements of the form $s'_J \otimes s'_j$, where $j \geq j_1$.

Let $W  \subset V_{g,l}$ be
the subspace generated by the elements of the form $s_I \otimes
s_j$, where $j \not \in I$ and $j \geq i_1$.

\begin{REM}
\label{W_i} The map $F_l :V_{g} \rightarrow \Lambda^{l-1}
\otimes B_g$ is injective if and only if ${F_l}_{|W}$ is
injective.
\end{REM}
\proof
The proof is completely analogous to the proof of \eqref{W}.
\qed

\begin{TEO}
\label{can}
If property $N_{p}$ holds for a canonical binary curve of genus $g
-1$  as in \eqref{can}, then the same property holds for a
canonical binary curve of genus $g $ as in \eqref{can} for
a generic choice of the parameters.
\end{TEO}
\proof

From the above discussion we know that the statement is equivalent to prove
injectivity of the Koszul map
$ F_{l}: V_{g} \rightarrow \Lambda^{l-1}H_g \otimes D_g$
for $l = g-2-p$, while by assumption we know injectivity of the map
$\tilde{F}_{l-1}: V_{g-1,r} \rightarrow \Lambda^{l-2}H_{g-1,r} \otimes D_{g-1,r}.$

We first project from $P_{g}$.
By Remark \eqref{W_i} we can reduce to prove injectivity of $F_{l}$ restricted the subspace $W$ generated by such $s_I \otimes s_j$ with $j \not \in I, j >i_1$. Note that as before we can decompose $W$ as  $W:=X_{g} \oplus Y_{g}$, where $X_{g}$ is the intersection with $W$ of the image of the map $\wedge s_{g} \otimes J_{g}$ in diagram \eqref{diagram1can} and $Y_{g}$ is the subspace of $W$ generated by such $s_I \otimes s_j$ with $g \not \in I$ and $j \not \in I, j >i_1$: $$X_{g} = \langle s_{g} \wedge  s_{J}  \otimes s_{j} \ | \   j \not \in J, j > j_1\rangle,  \ Y_{g} = \langle s_I \otimes s_{j} \ | \ j ,g \not \in I, j>i_1 \rangle$$

If  $F_{l}(x_{g} + y_{g}) =0$, where $x_{g} \in X_{g}$, $y_{g} \in Y_{g}$, then  $0=\pi_{g} \circ F_{l}(x_{g} + y_{g}) = \pi_{g} \circ F_{l}(x_{g})= (\wedge s_{g} \otimes L_{g}) \circ  \tilde{ F}_{l-1}(x_{g})$. Hence $x_{g}=0$, since by induction $\tilde{ F}_{l-1}$ is injective.
 So we have  reduced to prove injectivity of $F_{l}$ restricted to $Y_{g}$.

Repeat the procedure, i.e. project from the points $P_r$, $r=g-1...l$.  In this way we can reduce to prove injectivity for the  restriction of the map $F_{l}$ to the subspace $Y_{l}$ of $W$ generated by the elements of the form $s_I \otimes  s_j$ where $l,...,g,j \not \in I, j>i_1$. Observe that since $|I| = l$, we have  $Y_{l} = 0$, so $F_{l}$ is injective and the theorem is proved.

\qed

\begin{REM}
Notice that, by the theorem of Green and Lazarsfeld (\cite{gr}), if $p>g-[\frac{g}{2}]-2$, condition $N_p$ does not hold for any curve $\tilde{C}$ of genus $g-1$.
\end{REM}

\begin{COR}
If the Green conjecture is true for a canonical binary curve of
genus $g = 2k-1$  as in \eqref{can}, then it is true for a
canonical binary curve of genus $g = 2k$ as in \eqref{can}
for a generic choice of the parameters.
\end{COR}
\proof The conjecture for $g = 2k$ says that $K_{k,1}(C,K_C) = 0$,
or analogously that property $N_{k-2}$ holds for $C$ embedded with
$K_C$. By assumption we know that
$K_{k-1,1}(\tilde{C},K_{\tilde{C}}) = 0$, namely that property
$N_{k-2}$ holds for $\tilde{C}$ embedded with $K_{\tilde{C}} $, so
the thesis immediately follows from \eqref{can}
. \qed

With maple  (http://www-dimat.unipv.it/~frediani/greenfinal.tar.gz) one verifies the conjecture for $g = 5,7, 9,11$, so one can prove with the same method that if   $g \geq 3$, then propery $N_0$ holds (see also \cite{ccm} section 2), if $g \geq 5$, then propery $N_1$ holds, if $g \geq 7$, then propery $N_2$ holds, and if $g \geq 9$, then property $N_3$ holds, and if $g \geq 11$, then property $N_4$ holds.

\end{document}